Любимому папе посвящаю

# Элементарный вывод формулы Ю. В. Нестеренко разложения в цепную дробь числа $2 \times \zeta(3)$
## С.Н. Гладковский

Автором предложен элементарный вывод формулы Ю. В. Нестеренко разложения в цепную дробь числа $2 \times \zeta(3)$, обладающая хорошей скоростью сходимости.
Ключевые слова: непрерывная дробь, цепная дробь, $\zeta(3)$.

Как известно, $2 \times \zeta(3)$ можно представить разными видами цепной дроби. Например, используя формулу для $\zeta(3)$, данную Апери (см. напр. в [1] ):

$$\left.\begin{array}{l} 2\zeta(3) = \dfrac{12}{A_0} \\[2mm] A_k = 34k^3 + 51k^2 + 27k + 5 - \dfrac{(k+1)^6}{A_{k+1}} \end{array}\right\} \qquad (1)$$

В работе 1996г. [2] Ю. В. Нестеренко отмечает, что после работы Апери "За прошедшие с тех пор почти 20 лет не было получено ничего принципиально нового." Несколько далее, раскрывая содержание, указанный автор пишет: " Ниже (§1) мы предлагаем доказательство теоремы 1, инициированное работой Л.А. Гутника [7], и затем (§2) доказываем следующую ниже теорему 2, где даётся новое разложение $\zeta(3)$ в непрерывную дробь". На самом деле эта теорема о разложении $2 \times \zeta(3)$ в непрерывную дробь (причина этого будет объяснена ниже):

" Теорема 2. Справедливо следующее разложение в непрерывную дробь

$$2\zeta(3) = 2 + \cfrac{1}{2 + \cfrac{2}{4 + \cfrac{1}{3 + \cfrac{4}{2 + \cfrac{2}{4 + \cfrac{6}{6 + \cfrac{4}{5 + \cfrac{9}{4 + \cfrac{6}{6 + \cfrac{12}{8 + \cfrac{9}{7 + \cfrac{16}{6 + \ldots}}}}}}}}}}}, \quad (2)$$

где числители $a_n$, $n \geq 2$, и знаменатели $b_n$, $n \geq 1$, определяются по формулам

$$\begin{array}{ll} b_{4k+1} = 2k + 2, & a_{4k+1} = k(k+1), \\ b_{4k+2} = 2k + 4, & a_{4k+2} = (k+1)(k+2), \\ b_{4k+3} = 2k + 3, & a_{4k+3} = (k+1)^2, \\ b_{4k+4} = 2k + 2, & a_{4k+4} = (k+2)^2. \text{"} \end{array}$$

В приведённой выше формулировке наглядны недостатки нынешней общепринятой системы записи цепной дроби:

– громоздкость записи;

– разная индексация для числителей и знаменателей (для числителей $a_n$ отсчёт индекса начинается с 2, а для знаменателей $b_n$ – с 1) вместо единой, например, как у вашего покорного слуги, с 0;

– затушевана периодичность структуры дроби.

Две последние причины вынудили сформулировать теорему для $2\zeta(3)$, так как в случае для $\zeta(3)$ создавали дополнительную громоздкость и нарушали стройность формулировок и для числителей, и для знаменателей.

В работе (мне доступна версия только на английском) Л.А. Гутника [3], доложенной в октябре 2009г., представлено доказательство *двух других* цепных дробей для $2\zeta(3)$.

Излагая результаты Апери и Нестеренко, Гутник сравнивает уже отредактированную им цепную дробь Нестеренко для $\zeta(3)$. И приходит к выводу: " It is easy to prove that numerator and denominator of Nesterenko fractions with subscript $4v$–2 are equal to the numerator and denominator of Apery fraction with subscript $v$, respectively."

Собственно это замечание и побудило автора настоящей заметки вывести цепную дробь Нестеренко из так называемой цепной дроби Апери в виде (1) элементарными средствами.

Чтобы избежать путаницы в дальнейшем изложении, приведу результат Нестеренко в привычном для настоящей заметки виде:

$$\left.\begin{array}{l} 2\zeta(3) = 2 + \dfrac{1}{N_0} \\[2mm] N_k = 2k + 2 + \dfrac{(k+1)(k+2)}{2k+4 + \dfrac{(k+1)^2}{2k+3 + \dfrac{(k+2)^2}{2k+2 + \dfrac{(k+1)(k+2)}{N_{k+1}}}}} \end{array}\right\} \qquad (3)$$

Перепишем (1) в виде:

$$\left.\begin{array}{l} \zeta(3) = \dfrac{6}{5 - \dfrac{1^6}{A_0}} \\[3mm] A_k = 34(k+1)^3 + 51(k+1)^2 + 27(k+1) + 5 - \dfrac{(k+2)^6}{A_{k+1}} \end{array}\right\} \qquad (4)$$

или

$$A_k = 34k^3 + 153k^2 + 231k + 117 - \frac{(k+2)^6}{A_{k+1}} \qquad (5)$$

или

$$A_k = 5(k+1)^3 + 29k^3 + 138k^2 + 216k + 112 - \frac{(k+2)^6}{A_{k+1}}. \qquad (6)$$

Пусть $W_k = T_k - 5(k+1)^3$, тогда (6) примет вид:

$$W_k = 29k^3 + 138k^2 + 216k + 112 - \frac{(k+2)^6}{W_{k+1} + 5(k+2)^3}$$

или

$$W_k = 29k^3 + 138k^2 + 216k + 112 - \frac{(k+2)^6}{W_{k+1} + (5k^2 + 20k + 20)(k+2)}. \qquad (7)$$

Запишем (7) в виде:

$$W_k = 29k^3 + 138k^2 + 216k + 112 - \frac{(k+2)^5}{6 * \dfrac{W_{k+1}}{6(k+2)} + (5k^2 + 20k + 20)}$$

Положим $U_k = \dfrac{W_k}{6(k+1)}$, тогда имеем:

$$U_k = \frac{29k^3 + 138k^2 + 216k + 112}{6(k+1)} + \frac{1}{6(k+1)} \frac{-(k+2)^5}{(5k^2 + 20k + 20) + 6U_{k+1}} \qquad (8)$$

Нетрудно убедиться, что (8) можно представить в виде четырёхступенчатой цепной

дроби

$$U_k = (2k+3)(2k+4) + \cfrac{(k+2)^3}{(k+1) + \cfrac{(k+1)}{4 + \cfrac{1}{1 + \cfrac{1}{\cfrac{U_{k+1}}{(k+2)^2}}}}}. \qquad (9)$$

Полагая в (9) $P_k = \dfrac{U_k}{(k+1)^2}$, получаем

$$P_k = \cfrac{(2k+3)(2k+4)}{(k+1)^2} + \cfrac{\cfrac{(k+2)^3}{(k+1)^3}}{1 + \cfrac{1}{4 + \cfrac{1}{1 + \cfrac{1}{P_{k+1}}}}}, \qquad (10)$$

а в (10) $Q_k = 1 + \dfrac{1}{4 + \cfrac{1}{1 + \cfrac{1}{P_k}}}$, имеем:

$$Q_k = 1 + \cfrac{1}{4 + \cfrac{1}{1 + \cfrac{1}{\cfrac{(2k+3)(2k+4)}{(k+1)^2} + \cfrac{(k+2)^3}{(k+1)^3}}{Q_{k+1}}}} \qquad (11)$$

или

$$\left.\begin{array}{l} \zeta(3) = Q_0 \\[2mm] Q_k = 1 + \cfrac{1}{4 + \cfrac{1}{1 + \cfrac{(k+1)^3}{(k+1)(2k+3)(2k+4) + \cfrac{(k+2)^3}{Q_{k+1}}}}} \end{array}\right\}, \qquad (12)$$

то есть формула автора настоящей заметки, приведённая без доказательства в [4], (потому как пришёл к этой формуле, читая [1] и [3]).

Запишем (12) в виде

$$\left.\begin{array}{l} 2\zeta(3) = 2Z_0 \\[2mm] Z_k = 1 + \cfrac{1}{4 + \cfrac{2}{2 + \cfrac{2(k+1)^3}{2(k+1)(k+2)(2k+3) + \cfrac{2(k+2)^3}{2Z_{k+1}}}}} \end{array}\right\}. \qquad (13)$$

Пусть $2Z_k = H_k$, тогда (13) примет вид

$$2\zeta(3) = H_0$$

$$H_k = 2 + \cfrac{1}{2 + \cfrac{1}{2 + \cfrac{(k+1)^3}{(k+1)(k+2)(2k+3) + \cfrac{(k+2)^3}{H_{k+1}}}}} \qquad (14)$$

Положим $H_k = 2 + \dfrac{1}{G_k}$, имеем:

$$G_k = 2 + \cfrac{1}{2 + \cfrac{(k+1)^3}{(k+1)(k+2)(2k+3) + \cfrac{(k+2)^3}{2 + \cfrac{1}{G_{k+1}}}}} \qquad (15)$$

или

$$G_k = 2 + \cfrac{k+2}{2k+4 + \cfrac{(k+1)^3}{(k+1)(2k+3) + \cfrac{(k+1)(k+2)^2}{2k+2 + \cfrac{k+1}{G_{k+1}}}}} \qquad (16)$$

или

$$2\zeta(3) = 2 + \frac{1}{G_0}$$

$$G_k = 2 + \cfrac{k+2}{2k+4 + \cfrac{(k+1)^2}{2k+3 + \cfrac{(k+2)^2}{2k+2 + \cfrac{k+1}{G_{k+1}}}}} \qquad (17)$$

Полагая $N_k = (k+1)G_k$ получим:

$$2\zeta(3) = 2 + \frac{1}{N_0}$$

$$N_k = 2k+2 + \cfrac{(k+1)(k+2)}{2k+4 + \cfrac{(k+1)^2}{2k+3 + \cfrac{(k+2)^2}{2k+2 + \cfrac{(k+1)(k+2)}{N_{k+1}}}}}$$

что и требовалось.

\* \* \* \* \*

## Литература

S.N. Gladkovskii

An elementary derivation of the formula of Yu. V. Nesterenko expansion in continued fraction of a number $2\zeta(3)$.

The author proposed an elementary derivation of the formula of Yu. V. Nesterenko expansion in continued fraction of a number $2\zeta(3)$.

Гладковский Сергей Николаевич

E-mail: journaly2010@bk.ru